\documentclass[12pt]{article}
\usepackage{amssymb, latexsym, amsmath, amsthm}
\newtheorem{corollary}{Corollary}
\newtheorem{theorem}{Theorem}

\newtheorem{lemma}{Lemma}

\newtheorem{claim}{Claim}

\begin{document}
\title{{\bf A Note on Ramsey Numbers for Books}}
\author{V. Nikiforov and C. C. Rousseau \\ Department of
Mathematical
Sciences \\ The University of Memphis \\ 373 Dunn Hall \\
Memphis, Tennessee 38152-3240}
\date{}
\maketitle
\begin{abstract}
The {\em book with} $n$ {\em pages} $B_n$ is the graph consisting
of $n$ triangles sharing an edge.  The {\em book Ramsey
number} $r(B_m,B_n)$ is the smallest integer $r$ such
that either $B_m \subset G$ or $B_n \subset \overline{G}$
for every graph $G$ of order $r$. We prove that there exists a
positive constant $c$ such that $r(B_m, B_n) = 2n+3$ for all $n
\ge c \,m$.
\end{abstract}

\setlength{\baselineskip}{.25in} \section{Introduction} The graph
$B_n = K_1 + K_{1,n}$, consisting of $n$ triangles sharing a
common edge, is known as the {\em book} with $n$ {\em pages}. The
book Ramsey number $r(B_m, B_n)$ is the smallest integer $r$ such
that either $B_m \subset G$ or $B_n \subset \overline{G}$ for
every graph $G$ of order $r$.  The
study of Ramsey numbers for books was initiated in \cite{RS:78}
and continued in \cite{FRS:82}.   The
following results are known.
\begin{theorem}[Rousseau, Sheehan] For all $n > 1$, $r(B_1, B_n)
= 2n+3$.
\end{theorem}
\begin{theorem}[Parsons, Rousseau, Sheehan] If $2(m+n+1) >
(n-m)^3/3$
then $r(B_m, B_n) \le 2(m+n+1)$. More generally, \begin{equation}
\label{eq:par} r(B_m, B_n) \le m+n+2 + \left\lfloor \tfrac{2}{3}
\sqrt{3(m^2+mn+n^2)} \right\rfloor. \end{equation} If $4n+1$ is a
prime power, then $r(B_n, B_n) = 4n+2$.  If $m \equiv 0 \pmod{3}$
then $r(B_m,B_{m+2}) \le 4m+5$.
\end{theorem}
The more general upper bound (\ref{eq:par}) was noted by Parsons
in \cite{Par:STGT78}.  In looking for cases where
equality holds in (\ref{eq:par}) or in other cases covered by
Theorem 1, it is natural to consider the class of strongly
regular graphs. A $(v,k,\lambda,\mu)$ {\em strongly regular
graph} (SRG) is a graph with $v$ vertices that is regular of
degree $k$ in which any two distinct vertices have $\lambda$
common neighbors if they are adjacent and $\mu$ common neighbors
if they are nonadjacent. Thus if a $(v,k,\lambda,\mu)$ graph
exists then
\[
r(B_{\lambda + 1}, B_{v - 2k + \mu - 1}) \ge v+1.
\]
Inspection of a table known strongly regular graphs
\cite{Brou:HCD96} yields a number of exact values for book Ramsey
numbers.
\begin{corollary} In addition those cases where $4n+1$
is a prime power and $r(B_n,B_n) = 4n+2 \;
(n=1,2,3,4,6,\ldots,69)$, Theorem 2 gives the following exact
values for $r(B_m, B_n)$ in which the lower bound comes from a
strongly regular graph of order at most $280$.
\begin{center}
\begin{tabular}{||c|c|c||} \hline\hline $(m,n)$ & $r(B_m,B_n)$ &
$(v,k,\lambda,\mu)$  \\ \hline\hline
(2,5) & 16 & (15,6,1,3)  \\
(3,5) & 17 & (16,6,2,2) \\
(4,6) & 22 & (21,10,3,6) \\
(7,10) & 36 & (35,16,6,8)  \\
(11,11) & 46 & (45,22,10,11) \\
(14,17) & 64 & (63,30,13,15) \\
(23,26) & 100 & (99,48,22,24) \\
(22,37) & 120 & (119,54,21,27) \\
(29,38) & 136 & (135,64,28,32) \\
(34,37) & 144 & (143,70,33,35) \\
(47,50) & 196 & (195,96,46,48) \\
(46,58) & 210 & (209,100,45,50) \\
(56,56) & 226 & (225,112,55,56) \\
(38,82) & 244 & (243,110,37,60) \\
(62,65) & 256 & (255,126,61,63) \\
(69,71) & 281 & (280,135,70,60) \\ \hline\hline
\end{tabular}
\end{center}
\end{corollary}
The starting point for this paper is Theorem 1 together with the
following pair of results from \cite{FRS:82}.
\begin{theorem}[Faudree, Rousseau, Sheehan]
\[
r(B_2, B_n) \le \begin{cases} 2n+6, & 2 \le n \le 11, \\
                              2n+5, & 12 \le n \le 22, \\
                              2n+4, & 23 \le n \le 37, \\
                              2n+3, & n \ge 38.
                              \end{cases}
\]
\end{theorem}
\begin{theorem}[Faudree, Rousseau, Sheehan] If $m > 1$ and
\[
n \ge (m-1)(16m^3 + 16m^2 - 24m - 10) + 1, \] then $r(B_m, B_n) =
2n+3$.
\end{theorem}
From these results, we see that for each $m$ there exists a
smallest positive integer $f(m)$ such that  $r(B_m, B_n) = 2n+3$
for all $n \ge f(m)$.  Moreover $f(1) = 2$ and $f(2) \le 38$. Our
main purpose here is to prove the following strengthening of
Theorem 4.
\begin{theorem} There exists a positive constant $c$ such that
$r(B_m, B_n) = 2n+3$ for all $n \ge cm$.
\end{theorem}
\section{Proofs}
For standard terminology and notation, see \cite{Bol:MGT98}. For
$v \in V(G)$  we denote the neighborhood of
$v$ by $N(v)$ and the degree of $v$ by $\deg(v)$.
If needed, we shall use a subscript to identify the graph in
question; for example, $N_G(v)$ denotes the neighborhood of $v$
in  $G$. Given two
disjoint sets $U,W \subset V(G)$, let
$e(U,W) = |\{uw \in E(G)| \; u \in U, \, w \in W\}|$.
 The subgraph of $G$ induced by $X \subset V(G)$ will be
denoted by $G[X]$. Given graphs $G$ and $H$, let $M_G(H)$ denote
the number of induced subgraphs of $G$ that are isomorphic to
$H$. The number of pages in the largest book contained in $G$
will be called the {\em book size} of $G$ and this will be
denoted by $bs(G)$.  It is convenient to identify the graph
and its complement in terms of edge colorings of a complete
graph.  In this framework, $r(B_m, B_n)$ is the smallest $r$ such
that in every $(R,B)$ = (red, blue) coloring of $E(K_r)$, either
$bs(R) \ge m$ or $bs(B) \ge n$.  Theorem 5 clearly follows from
the following fact.
\begin{theorem} Suppose $m$ and $n$ are positive integers
satisfying $n \ge 10^6 m$.  If $(R,B)$ is any two-coloring of
$E(K_n)$ then either $bs(R) > m$ or $bs(B) \ge n/2 - 2$.
\end{theorem}
In view of the case $R = K(\lfloor n/2 \rfloor, \lceil n/2
\rceil)$ the conclusion $bs(B) \ge n/2 - 2$ is best possible. The
proof of Theorem 6 uses the following counting result.
\begin{lemma} Let $G$ be a graph with $p$
vertices and $q$ edges that satisfies $bs(G) \le m$.  For $0 <
\lambda < 1$ suppose that $\delta(G) \ge \lambda p$ so $q \ge
\lambda p^2/2$. If $p> 5(2 \lambda + 1)/\lambda^2$ then
\[
M_G(C_4) > \left( \frac{\lambda^3 p^2}{5} - \frac{m^2}{2}
\right)q.
\]
\end{lemma}
\begin{proof}
For distinct vertices $u,v \in V(G)$ let $c(u,v) = \left|N_G(u)
\cap
N_G(v)\right|$.  Then
\begin{align*}
\sum_{\{u,v\}} \binom{c(u,v)}{2} & = 2 M_G(C_4) + 6M_G(K_4) + 2
M_G(B_2) \quad \text{and}
\\[.1in]
\sum_{uv \in E} \binom{c(u,v)}{2} & = 6M_G(K_4) + M_G(B_2),
\end{align*}
from which we get \begin{equation} \label{eq:C4} M_G(C_4) =
\frac{1}{2} \sum_{\{u,v\}} \binom{c(u,v)}{2} - \sum_{uv \in E}
\binom{c(u,v)}{2} + 3M_G(K_4).
\end{equation}
Since $bs(G) \le m$,
\[
\sum_{uv \in E} \binom{c(u,v)}{2} \le q \binom{m}{2} <
\frac{qm^2}{2}.
\]
Note that
\[
\sum_{\{u,v\}} c(u,v) = \sum_{v \in V(G)} \binom{\deg_G(v)}{2}
\ge p \binom{2q/p}{2} = q(2q/p-1) \ge q(\lambda p-1) := x,
\]
so by convexity,
\[
\sum_{\{u,v\}} \binom{c(u,v)}{2} \ge \binom{p}{2}
\binom{x/\binom{p}{2}}{2} = \frac{x}{2} \left(
\frac{x}{\binom{p}{2}} - 1 \right). \] Since \[
\frac{x}{\binom{p}{2}}
> \frac{2x}{p^2} = \frac{2q(\lambda p -1)}{p^2} \ge
\lambda(\lambda p - 1), \] we
have \[ \sum_{\{u,v\}} \binom{c(u,v)}{2} > \frac{q(\lambda p -
1)(\lambda(\lambda p - 1) - 1)}{2} > \frac{2 \lambda^3 p^2 q}{5}.
\] {\em Note.} The last inequality is clear if $\lambda^2 p - 10
\lambda  - 5 \ge 0$,
and hence it holds since we have required $p \ge 5(2\lambda +
1)/\lambda^2$. In view of (\ref{eq:C4}) we have
\[
M_G(C_4) > \left(\frac{\lambda^3 p^2}{5} - \frac{m^2}{2}
\right)q,
\]
as claimed.
\end{proof}
\begin{proof}[Proof of Theorem 6]  Suppose $n \ge 10^6 m$ and
that $(R,B)$ is a two-coloring of $E(K_n)$ such that $bs(R) \le
m$.  We shall prove that $bs(B) \ge n/2 - 2$.
Let ${\cal H} = C_4 \cup K_1$.
\begin{claim}
If $bs(B) \le n/2 - 2$ then $M_R({\cal H}) \le 4m M_R(C_4)$.
\end{claim}
{\em Note.}  The hypothesis $bs(G) \le n/2 - 2$ rather than, as
one might naturally expect, $bs(G) < n/2 - 2$, is made for 
convenience.
\begin{proof}
Suppose $M_R({\cal H}) > 4m M_R(C_4)$. Then there exists an
induced $C_4 = (u,v,w,z)$ such that
\[
|N_{B}(u) \cap N_{B}(v) \cap N_{B}(w) \cap N_{B}(z)| \ge 4m+1.
\]
Since $bs(B) \le n/2 - 2$ we have \[ |N_{B}(u) \cap N_{B}(w)| \le
n/2 - 2 \qquad \text{and} \qquad |N_{B}(v) \cap N_{B}(z)| \le n/2
- 2.
\]
It then follows that there are at least $4m+1$ vertices outside
of $\{u,v,w,z\}$ that are adjacent in $R$ to at least one of
$u,w$
and at least one of $v,z$. This gives $m+1$ or more red triangles
on at least one of the four edges $uv, vw, wz, zu$, and thus the
desired contradiction.
\end{proof}
It is known that for any graph $G$ of order $n$,
\[
M_G(C_4) \le \binom{\lfloor n/2 \rfloor}{2} \binom{\lceil n/2
\rceil}{2} < \frac{n^4}{64}.
\]
See \cite{BNT:86} for a proof of the more general result \[
M_G(K_{m,m}) \le \binom{\lfloor n/2 \rfloor}{m} \binom{\lceil n/2
\rceil}{m}.
\]
 Hence by Claim 1,
\[
M_R({\cal H}) < \frac{m n^4}{16}
\]
or else $bs(B) > n/2 - 2$.
\begin{claim}
If $bs(B) \le n/2 - 2$ then $R$ has at most $n/20$
vertices of degree $9n/20$ or less.
\end{claim}
\begin{proof} Let $v$ be any vertex of degree $9n/20$ or less in
$R$ and let $X = N_B(v)$. Then $B[X]$ has maximum degree at most
$n/2 - 2$ so $G = G(v) = R[X]$ has minimum degree $\delta$
satisfying
\[
\delta \geq |X| - 1 - \frac{n}{2} + 2 \ge n - 1 - \frac{9n}{20} -
1 - \frac{n}{2} + 2 =  \frac{n}{20}
\]
and (since $|X|+1 \ge 11n/20$) \[ \delta \ge |X|+1 - \frac{n}{2}
\ge |X|+1 - \frac{1}{2} \left(\frac{20(|X|+1)}{11} \right) =
\frac{|X|+1}{11}.
\] Let us check that Lemma 1 applies to $G$.  Take $\lambda =
1/11$ and $p = |X| \ge 11n/20$. Then $p > 5(2 \lambda +
1)/\lambda^2$ holds provided $n
\ge 1302$. This is certainly the case since $n \ge 10^6 m$. Using
$m \le n/10^6$, Lemma 1 gives
\[
M_G(C_4) > \left( \frac{1}{5} \cdot \frac{1}{11^3}
\left(\frac{11n}{20}\right)^2 - \frac{1}{2}\left( \frac{n}{10^6}
\right)^2 \right) \frac{1}{2} \left( \frac{11n}{20} \right)
\left(
\frac{n}{20} \right) \approx \frac{n^4}{640,000}.
\]
Suppose more than $n/20$ vertices in $R$ have degree
$9n/20$ or less.  Then \[ \frac{m n^4}{16} > M_R({\cal H}) =
\sum_{v} M_{G(v)}(C_4)  > \sum_{\deg(v) \le 9n/20} M_{G(v)}(C_4)
> \frac{n}{20} \cdot \frac{n^4}{640,000},
\]
so $n < 8 \cdot 10^5m$, a contradiction.
\end{proof}
Let $S = \{v| \; \deg_R(v) > 9n/20\}$. From Lemma 3 we know that
$|S| > 19n/20$, so the minimum degree of $R[S]$ satisfies
\[
\delta \geq \frac{9n}{20} - (n - |S|) > \frac{2n}{5} \geq
\frac{2|S|}{5}.
\]
Now we use the following result of Andr\'{a}sfai, Erd\H{o}s and
S\'{o}s \cite{AES:74}.
\begin{theorem}[Andr\'{a}sfai, Erd\H{o}s, S\'{o}s]
Suppose $r \ge 3$. For any graph $G$ of order $n$,
at most two of the following properties can hold:
\[
\text{(i)} \;\;  K_r \not\subseteq G, \qquad \text{(ii)} \;\;
\delta(G)
> \frac{3r-7}{3r-4} \, n, \qquad \text{(iii)} \;\; \chi(G) \ge r.
\]
\end{theorem}
{\em Note.} In particular, a triangle-free graph $G$ with
$\delta(G) > 2|V(G)|/5$ is bipartite.

Now we are prepared to complete the proof of Theorem 6. It is
easy to see that $R[S]$ has no triangle. If
$T=\{u,v,w\}$ is a triangle in $R[S]$ and $U$ is the set of $n-3$
vertices outside $T$, then
\[
3(9n/20 - 2) < e_R(T,U) \le 3(m-1) \cdot 2 + (n-3(m-1)) = n +
3(m-1),
\]
or $7n/20 < 3m+3$, which is false.  Hence $R[S]$ is
bipartite by Theorem 7. Let $S_1$ and $S_2$ denote the two color
classes of $R[S]$.  Put $v \in T_1$ if $v$ is adjacent in $B$ to
every vertex of $S_1$.  Then for the remaining vertices put
$v \in T_2$ if $v$ is adjacent in $B$
to every vertex of $S_2$. Let $W_1 = S_1 \cup T_1, W_2 = S_2 \cup
T_2$, and let $X$ denote the set of vertices in neither $W_1$ nor
$W_2$.  If $X = \varnothing$ then we may assume that $|W_1| \ge
n/2$. In this case it is clear that $bs(B) \ge n/2 - 2$.

We are left to consider the case $X \ne \varnothing$. For
$u \in S$ let $Z(u) = N_B(u) \cap X$. For distinct vertices $u,v
\in S_1$, consideration of the blue book on $uv$ shows that
\begin{align*}
bs(B) & \ge |S_1| - 2 + |T_1| + |Z(u) \cap Z(v)| \\ & \ge |S_1| -
2 + |T_1| + |Z(u)| + |Z(v)| -|X|.
\end{align*}
Summing over all pairs $u,v \in S_1$ and computing the average,
we find \begin{align*} bs(B) & \ge |S_1| - 2 + |T_1| +
\frac{2(|S_1||X| - e_R(S_1,X))}{|S_1|} - |X| \\[.1in]
& = |S_1| + |T_1| + |X| - 2 - \frac{2e_R(S_1,X)}{|S_1|}.
\end{align*}
Similarly,
\[
bs(B) \ge |S_2| + |T_2| + |X| - 2 - \frac{2e_R(S_2,X)}{|S_2|}.
\]
Note that $|S_1| < n/2$ or else we are done at the outset;
similarly $|S_2| < n/2$.
Hence
\[ |S_1| = |S| - |S_2| > \frac{19n}{20} - \frac{n}{2} =
\frac{9n}{20},
\]
and likewise $|S_2| > 9n/20$. Consequently \begin{align*} bs(B) &
> |S_1| + |T_1| + |X| - 2 - \frac{40e_R(S_1, X)}{9n}, \\ bs(B) &
>
|S_2| + |T_2| + |X| - 2 - \frac{40e_R(S_2,X)}{9n}.
\end{align*}
Addition then gives
\[
2 \, bs(B) > n-4 + |X| - \frac{40e_R(S,X)}{9n}.
\]
Hence $e_R(S,X) > 9n|X|/40$ or else the proof is complete.

Thus we assume $e_R(S,X) > 9n|X|/40$ and now seek a companion
bound on $e_R(S,X)$.   For each $x \in X$ there is at least one
$v \in S_1$ such that $xv \in R$.
 Since $|N_R(v) \cap
S_2| \geq 2|S|/5$, consideration of the red book on $xv$ shows
that
\begin{align*}
bs(R) & \ge |N_R(x) \cap N_R(v) \cap S_2| \\
     & = |N_R(x) \cap S_2| + |N_R(v) \cap S_2| - |S_2| \\
     & \ge |N_R(x) \cap S_2| + \frac{2|S|}{5} - |S_2|.
     \end{align*}
Taking the average over $x \in X$, we obtain
\[
bs(R) \ge \frac{e_R(S_2,X)}{|X|} + \frac{2|S|}{5} - |S_2|.
\]
In exactly the same way,
\[
bs(R) \ge \frac{e_R(S_1,X)}{|X|} + \frac{2|S|}{5} - |S_1|.
\]
Hence
\[
2m \ge 2bs(R) \ge \frac{e_R(S,X)}{|X|} - \frac{|S|}{5}.
\]
Thus \[ e_R(S,X) \ge 2m|X| + \frac{|S||X|}{5}. \]
From the two bounds for $e_R(S,X)$, we obtain
\[
\frac{9n|X|}{40} < e_R(S,X) \le \frac{|S| |X|}{5} + 2m|X| <
\frac{n|X|}{5} + 2m|X|.
\]
By assumption $|X| > 0$, so
\[
\frac{9n}{40} < \frac{n}{5} + 2m,
\]
which is false. \end{proof}
\section{Concluding Remarks} The determination of the best
constant $c$ in Theorem 5 is open, as are other basic problems on
book Ramsey numbers stated in \cite{FRS:82}.  In particular, it
is unknown whether or not there exists a constant $C$ such that
$r(B_m, B_n) \le 2(m+n)+C$ for all $m,n$.

\vspace*{.3in}

\noindent
V. Nikiforov  \\
\texttt{vlado\_nikiforov@hotmail.com} \\[.2in]
C. C. Rousseau \\
\texttt{ccrousse@memphis.edu}

\begin{thebibliography}{9}
\bibitem{AES:74} B. Andra\'{a}sfai, P. Erd\H{o}s, and V. T.
S\'{o}s,
On the connection between chromatic number, maximal clique and
minimal degree of a graph, Discrete Math. 8 (1974), 205--218.
\bibitem{Bol:MGT98} B. Bollob\'{a}s, {\em Modern Graph
Theory}, Springer-Verlag, New York, 1998.
\bibitem{BNT:86} B. Bollob\'{a}s, C. Nara, and S. Tachibana, The
maximal number of induced complete bipartite graphs, Discrete
Math. 62 (1986), 271--275.
\bibitem{Brou:HCD96} A. E. Brouwer, ``Strongly regular graphs,''
Handbook of Combinatorial Designs, C. J. Colbourn and J. H.
Dinitz, (Editors), CRC Press, Boca Raton, 1996, pp. 667--685..
\bibitem{FRS:82} R. J. Faudree, C. C. Rousseau, and J. Sheehan,
Strongly regular graphs and finite Ramsey theory, Linear Algebra
Appl.  46 (1982), 221--241. \bibitem{Par:STGT78} T. D. Parsons,
``Ramsey graph theory,'' Selected Topics in Graph Theory, L. W.
Beineke and R. J. Wilson, (Editors), Academic Press, London,
1978, pp. 361--384.
\bibitem{RS:78} C. C. Rousseau and J. Sheehan, On Ramsey numbers
for
books, J. Graph Theory  2  (1978), 77--87.
\end{thebibliography}
\end{document}